\documentclass[letterpaper,10pt]{amsart}

\newcommand{\R}{\mathbb{R}}
\newcommand{\Z}{\mathbb{Z}}

\newtheorem{lm}{Lemma}[section]
\newtheorem{thq}[lm]{Theorem}
\theoremstyle{remark}
\newtheorem{rem}[lm]{Remark}

\DeclareMathOperator{\im}{Im}

\title[Automorphism of A Free Metabelian Group]{An Automorphism of A Free Metabelian Group Without Fixed Points}
\author{Martin Kassabov}

\address{
Department of Mathematics \\
University of Alberta\\
632 Central Academic Building\\
Edmonton, Alberta, T6G 2G1\\
Canada}

\email{kassabov@aya.yale.edu}
\urladdr{http://www.math.ualberta.ca/\~{}mkassabov}


\begin{document}
\setlength{\parindent}{0pt}


\begin{abstract}
We construct an example of an IA-automorphism of the free
metabelian group of rank $n\geq 3$ without
nontrivial fixed points. That gives a negative answer to
the question raised by Shpilrain in \cite{Sh}.
By a result of Bachmuth~\cite{Ba1}, such an automorphism
does not exist if the rank is equal to $2$.
\end{abstract}


\maketitle

\begin{section}{Introduction}
\label{sec:intro}

The study of free metabelian groups began in the 30-es with
the work of Magnus~\cite{Ma}. Later Bachmuth wrote a series of
papers~\cite{Ba1,Ba2,Ba3} devoted
to automorphisms of the free metabelian groups.
In~\cite{Sh} Shpilrain studied the fixed points of
the endomorphisms of free metabelian groups. There, he raised the question
whether any IA-automorphism (i.e., an automorphism which acts
trivially on the abelianization) of the a metabelian group has a
nontrivial fixed point. In this paper we construct explicitly
an IA-automorphism without fixed points when the number of
generators is at least $3$, thus giving negative answer to
Shpilrain's question. Such an automorphism does not exist in
the two-generator case by a result of Bachmuth.

In the next section we have summarized several known results
about free metabelian groups. Section~\ref{sec:auto}
describes an automorphism without fixed points.

The author thanks his adviser E. Zelmanov for useful discussions
during the preparation of the manuscript, and the referee for
suggesting significant refinements of remark~\ref{lie}.

\end{section}


\begin{section}{Free Metabelian Group}
\label{sec:gr}

Let $F_n$ be the free group on the generators $g_1,\dots, g_n$.
Let $F'_n = [F_n,F_n]$ be its commutator subgroup, and
let $F''_n = [F'_n,F'_n]$
be the second commutator subgroup. The group
$M_n = F_n /F''_n$ is called the free metabelian group of
rank $n$. We will denote by $g_i$ the generators of $M_n$.

Let $\phi$ be the two dimensional representation of the group
$M_n$ defined by
$$
\phi(g_i) = \left(\begin{array}{cc} s_i & t_i \\ 0 & 1 \end{array}
\right).
$$
The coefficients of these matrices lie in the ring
$\Z[s_i^{\pm 1},t_i]$ obtained
from $\Z$ by adjoining the commuting variables $s_i$ and $t_i$ together with
the inverses of $s_i$. This representation is known as the Magnus
representation. Magnus~\cite{Ma} showed that:

\begin{thq}
\label{gr:magnus}
The representation $\phi$ is faithful.
\end{thq}

The following lemma from~\cite{Ba1} describes the image of the
representation $\phi$.

\begin{lm}
\label{gr:image}
Let $S=s_1^{j_1} \dots s_n^{j_n}$, where $j_i$ are arbitrary integers.
Let $\gamma_i$ be polynomials in $s_i^{\pm 1}$. Then the matrix
$$
\left(\begin{array}{cc}
S & \sum \gamma_i t_i \\
0 & 1
\end{array} \right)
$$
lies in the image of $\phi$ if and only if $\gamma_i$ satisfy the identity
$$
\sum \gamma_i(1-s_i) = 1 - S.
$$
Also, every element in $\im(\phi)$ is a matrix of the above type.
\end{lm}

The following lemma (also from~\cite{Ba1}) gives a representation of the
IA-endomorphisms by matrices over the ring $\Z[s_i^{\pm 1}]$:

\begin{lm}
\label{gr:matrix}
Let $\alpha$ be an IA-endomorphism of the group $M_n$. There exist
polynomials $\alpha_{i,j}$ on $s_i^{\pm 1}$
such that
$$
\phi(\alpha(g)) = \bar \alpha (\phi(g)),
$$
where $\bar\alpha$ is the endomorphism of the coefficient ring
$\Z[s_i^{\pm 1}, t_i]$ of the representation $\phi$ defined by
$$
\begin{array}{r@{\,=\,}l}
\bar\alpha(t_i) & \sum \alpha_{i,j}t_j\\
\bar\alpha(s_i) & s_i.
\end{array}
$$
\end{lm}

\end{section}


\begin{section}{Construction of an Automorphism}
\label{sec:auto}

In this section we give an explicit example of an automorphism of
the group $M_n$, which has no fixed points.

\begin{lm}
Let $\alpha_n$ be an endomorphism of $M_n$ for $n\geq 3$
defined by
$$
\begin{array}{ll}
\alpha_n(g_i) = \big[[g_1,g_2]g_n,g_i \big]g_i &\mbox{if }i<n \\
\alpha_n(g_n) = [g_1,g_2]g_n. &\\
\end{array}
$$
Then $\alpha_n$ is an IA-automorphism of $M_n$.
Here, $[g,h]=ghg^{-1}h^{-1}$ denotes the commutator of the elements
$g$ and $h$.
\end{lm}
\begin{proof}
Let $\beta_1$ be the automorphism on $M_n$ which is the
conjugation by the element $g_n$, i.e., $\beta_1: g\to g_n g g_n^{-1}$; and let
$\beta_2$ be the automorphism defined by $\beta_2(g_i)=g_i$ for all
$i<n$ and $\beta_2(g_n) = [g_1,g_2]g_n$. The composition
$\beta_2 \circ \beta_1$ acts on the generators as follows:
$$
\begin{array}{ll}
\beta_2 \circ \beta_1(g_i) = \beta_2([g_n,g_i]g_i) = \big[[ g_1,g_2]g_n,g_i \big]g_i &\mbox{if }i<n \\
\beta_2 \circ \beta_1(g_n) = \beta_2(g_n) = [g_1,g_2]g_n &\mbox{for  }i=n. \\
\end{array}
$$
We have
$\alpha_n = \beta_2 \circ \beta_1$, and therefore $\alpha_n$ is an
automorphism. The automorphism $\alpha_n$ acts trivially on
the abelianization of $M_n$ because both $\beta_1$ and $\beta_2$ do.
\end{proof}

\begin{lm}
The automorphism $\alpha_n$ has no nontrivial fixed points in $M_n$.
\end{lm}
\begin{proof}
Let us give a direct proof using the Magnus representation.

Assume the contrary: let $g \not =1$ be a fixed point of the automorphism
 $\alpha_n$.
First, let us compute the action of the
endomorphism $\bar\alpha_n$ (defined in lemma~\ref{gr:matrix})
corresponding to the automorphism $\alpha_n$.

\begin{lm}
The endomorphism $\bar\alpha_n$ acts on $t_i$ as follows:
$$
\begin{array}{r@{\,=\,}l}
\bar \alpha (t_i) & s_n t_i + (1-s_i)\big((1-s_2) t_1 - (1-s_1) t_2 + t_n\big),  \mbox{ for } i \leq n \\
\bar \alpha (t_n) & (1-s_2) t_1 - (1-s_1) t_2 + t_n.
\end{array}
$$
\end{lm}
\begin{proof}
We have that
$$
\phi([g_1,g_2]) \! = \!
\left(\begin{array}{cc}
s_1 & t_1 \\
0 & 1
\end{array} \right) \!
\left(\begin{array}{cc}
s_2 & t_2 \\
0 & 1
\end{array} \right) \!
\left(\begin{array}{cc}
s_1^{-1} & -s_1^{-1}t_1 \\
0 & 1
\end{array} \right) \!
\left(\begin{array}{cc}
s_2^{-1} & -s_2^{-1}t_2 \\
0 & 1
\end{array} \right)
$$
$$
=
\left(\begin{array}{cc}
1 & (1-s_2) t_1 - (1-s_1) t_2 \\
0 & 1
\end{array} \right).
$$
Therefore,
$$
\phi(\alpha_n(g_n)) =
\left(\begin{array}{cc}
s_n & (1-s_2) t_1 - (1-s_1) t_2 + t_n \\
0 & 1
\end{array} \right)
$$
and for $i\leq n$
$$
\phi(\alpha_n(g_i)) =
\phi(\alpha_n(g_n))
\left(\begin{array}{cc}
s_i & t_i\\
0 & 1
\end{array} \right)
\phi(\alpha_n(g_n))^{-1}
=
$$
$$
=
\left(\begin{array}{cc}
s_i & s_n t_i + (1-s_i)\big((1-s_2) t_1 - (1-s_1) t_2 + t_n\big) \\
0 & 1
\end{array} \right).
$$
This, combined with the definition of $\bar \alpha_n$, proves the lemma.
\end{proof}

Let
$$
\phi(g)=
\left(\begin{array}{cc}
S & \sum \gamma_i t_i \\
0 & 1
\end{array} \right),
$$
where $\gamma_i$ are polynomials in $s_i^{\pm 1}$.

\begin{lm}
If $\phi(g) = \bar \alpha_n (\phi(g))$, then $\gamma_i=0$ for $3\leq i < n$, and
there exists $A \in \Z[s_i^{\pm 1}]$, such that
$$
\begin{array}{r@{\,=\,}l}
\gamma_1 & (1-s_2) s_n A, \\
\gamma_2 & -(1-s_1) s_n A, \\
\gamma_n & (1-s_n) A.
\end{array}
$$
\end{lm}
\begin{proof}
We have that
$$
\begin{array}{l}
\bar \alpha_n (\sum \gamma_i t_i) =
\sum \gamma_i \bar \alpha_n(t_i)= \\
=\sum_{i\leq n} \gamma_i \bigg(s_n t_i + (1-s_i)\big((1-s_2) t_1 - (1-s_1) t_2 + t_n\big) \bigg)
+
\\
\quad\quad\quad \quad\quad\quad\quad\quad\quad
+ \gamma_n\big((1- s_2) s_n t_1 - (1-s_1) s_n t_2 +  t_n\big).
\end{array}
$$
Equating the coefficients of $t_i$ in $\sum \gamma_i t_i$ and
$\bar \alpha_n(\sum \gamma_i t_i)$, we have:
$$
\begin{array}{@{\mbox{ coefficient at }}l@{:\quad\!\!}r@{\,\,=\,\,}l}
t_1 & \gamma_1 & s_n \gamma_1 + \sum_{j<n} (1-s_i)(1-s_2)\gamma_k + (1-s_2) s_n \gamma_n   \\
t_2 & \gamma_2 & s_n \gamma_2 - \sum_{j<n} (1-s_i)(1-s_1)\gamma_k - (1-s_1) s_n \gamma_n  \\
t_i & \gamma_i & s_n \gamma_i,\,\,  \mbox{ for } 3\leq i < n  \\
t_n & \gamma_n & \sum_{j < n} (1-s_j)\gamma_j + \gamma_n.
\end{array}
$$
The element $1-s_n$ is not a zero divisor in the ring $Z[s_i^{\pm 1}]$.
Therefore, we have $\gamma_i=0$, for $3\leq i <n$.
Thus, the system is reduced to
$$
\begin{array}{r@{\,=\,}l}
(1 - s_n) \gamma_1 & ( 1- s_2 ) \bigg(\big((1-s_1)\gamma_1 + (1-s_2)\gamma_2 \big)+ s_n\gamma_n \bigg) \\
(1 - s_n) \gamma_2 & -(1 - s_1) \bigg(\big((1-s_1)\gamma_1 + (1-s_2)\gamma_2 \big)+ s_n\gamma_n \bigg) \\
0 & (1 - s_1) \gamma_1 +(1- s_2) \gamma_2 .
\end{array}
$$
The expression in the middle brackets of the first two equations is zero. Therefore,
they are reduced to
$$
\begin{array}{r@{\,=\,}l}
(1 - s_n) \gamma_1 & ( 1- s_2 ) s_n\gamma_n \\
(1 - s_n) \gamma_2 & -(1 - s_1) s_n\gamma_n. \\
\end{array}
$$
The elements $1-s_i$ are prime in the ring $Z[s_i^{\pm 1}]$. Therefore,
$(1 - s_n)$ divides $\gamma_n$, $(1 - s_1)$ divides $\gamma_2$ and
$(1 - s_2)$ divides $\gamma_1$, i.e., there exist elements $A_i$ such that
$\gamma_1 =(1 - s_2) A_1$, $\gamma_2 = (1-s_1) A_2$ and
$\gamma_n = (1-s_n) A_n$. From the system above it follows that
$A_1= - A_2=s_n A_n$, which completes the proof.
\end{proof}

\begin{lm}
If $\phi(g) = \bar \alpha_n (\phi(g))$, then $\phi(g)=1$.
\end{lm}
\begin{proof}
The matrix
$\left(\begin{array}{cc}S & \sum \gamma_i t_i \\
0 & 1\end{array} \right)$
is in the image of $\phi$. By lemma~\ref{gr:image} we have that
$1-S=\sum \gamma_i (1-s_i)$. Applying the previous lemma, we have
$$
1-S = (1-s_n)^2 A.
$$
If $S\not=1$, then the point $(1,\dots,1)$ is a simple zero of
the left hand  side, but it is a double of the right hand side, a contradiction.
Therefore, the only solution to the above equation is $S=1$
and $A=0$, i.e., $\phi(g)=1$.
\end{proof}

This shows that any fixed point of the automorphism $\alpha_n$ lies in
the kernel of $\phi$. By lemma~\ref{gr:magnus} the representation
$\phi$ is faithful. Therefore, the automorphism $\alpha_n$ has no
nontrivial fixed points.
\end{proof}

Thus, we have shown:
\begin{thq}
There exists an automorphism of the free metabelian group on more
than $2$ generators which
acts trivially on the abelianization, and the only fixed point is the
identity.
\end{thq}

\begin{rem}
Every IA-automorphism of the group $M_2$ has a fixed
point. Bachmuth~\cite{Ba1} has shown that every such
automorphism is inner, and therefore, has a fixed point.
\end{rem}

\begin{rem}
\label{lie}
For the construction of the automorphism $\alpha_n$ we were motivated by the
Lie correspondence. It is well known that many problems regarding a Lie
group $G$ can be transferred to similar problems about its Lie algebra.

In general, these two problems are not equivalent, because the exponential map
is not a bijection. In the case of connected unipotent Lie groups, the
exponential map is a bijection, and one can use the Lie correspondence to relate
the question of existence of IA-automorphisms without fixed points to
the existence of derivations of their Lie algebras with trivial kernel.

We could not apply this correspondence directly to our problem because $M_n$ is not
a unipotent Lie group. One way to go around this problem is to
embed $M_n$ into a pro-unipotent Lie group
$\widetilde M_n$, where $\widetilde M_n$ is constructed as follows:
Let $A = \R \langle\langle x_1, \dots,x_n \rangle\rangle$
be the algebra of formal power series
of the non-commuting variables $x_i$ and let $A''$ be its second commutator subalgebra.
Finally, let $\widetilde M_n$ be the closure of the subgroup of the multiplicative group of
$1+ A/A''$ generated by the
images of $1+x_i$ together with the operations $g \to g^\alpha$ for
any real number $\alpha$, where the operation is defined as a formal power series.

The Lie algebra corresponding to the pro-unipotent group $\widetilde M_n$ is
the free metabelian Lie algebra $ML_n$ on $n$
generators, i.e., the quotient of the free Lie
algebra by its second commutator algebra.
Our first step was to construct a derivation $D_n$ of $ML_n$ with a trivial kernel ---
we took $D_n(x_i) = [x_i,x_n]$ for $i<n$ and $D_n(x_n) = [x_1,x_2]$
and proved that it has a trivial kernel using a well known basis of the free
metabelian algebra $ML_n$. The next step was to lift this derivation to
the automorphisms of $\widetilde M_n$ and $M_n$, i.e., to construct an automorphism
$\alpha_n$ of
$M_n$ such that its extension to $\widetilde M_n$ is almost the same as $\exp(D_n)$.

In the case of the free metabelian group,
the Lie correspondence is relatively weak ---
for
example there exist derivations of $ML_2$ without kernel but
they can not be used to construct automorphisms of the group $M_2$.

One way of showing that the automorphism $\alpha_n$ moves any nontrivial
element $g\in M_n$ is to pass to a suitable nilpotent quotient of $M_n$
and use the Lie correspondence there. Although such a proof would be
conceptually easier, the are many additional details that have to be proved, which would
make it more technical.  Therefore, we chose the more direct proof involving
the Magnus representation of the free metabelian group.
\end{rem}

\end{section}


\bibliographystyle{amsplain}
\bibliography{auto1}

\providecommand{\bysame}{\leavevmode\hbox to3em{\hrulefill}\thinspace}
\providecommand{\MR}{\relax\ifhmode\unskip\space\fi MR }
\providecommand{\MRhref}[2]{%
  \href{http://www.ams.org/mathscinet-getitem?mr=#1}{#2}
}
\providecommand{\href}[2]{#2}
\begin{thebibliography}{1}

\bibitem{Ba1}
S.~Bachmuth, \emph{Automorphisms of free metabelian groups}, Trans. Amer. Math.
  Soc. \textbf{118} (1965), 93--1104. \MR{31 \#4831}

\bibitem{Ba2}
S.~Bachmuth and H.~Y. Mochizuki, \emph{The nonfinite generation of ${\rm
  {a}ut}({G})$, ${G}$\ free metabelian of rank $3$}, Trans. Amer. Math. Soc.
  \textbf{270} (1982), no.~2, 693--700. \MR{83f:20026}

\bibitem{Ba3}
Seymour Bachmuth and Horace~Y. Mochizuki, \emph{${\rm {a}ut}({F})\to{\rm
  {a}ut}({F}/{F}'')$ is surjective for free group ${F}$ of rank $\geq 4$},
  Trans. Amer. Math. Soc. \textbf{292} (1985), no.~1, 81--101. \MR{87a:20032}

\bibitem{Ma}
Wilhelm Magnus, \emph{On a theorem of {M}arshall {H}all}, Ann. of Math.
  \textbf{40} (1939), 764--768. \MR{1,44b}

\bibitem{Sh}
Vladimir Shpilrain, \emph{Fixed points of endomorphisms of a free metabelian
  group}, Math. Proc. Cambridge Philos. Soc. \textbf{123} (1998), no.~1,
  75--83. \MR{98k:20056}

\end{thebibliography}


\end{document}